# JOINT DENSITY FOR THE LOCAL TIMES OF CONTINUOUS-TIME MARKOV CHAINS

By David Brydges,[1] Remco van der Hofstad[2]
and Wolfgang König[3]

*University of British Columbia, Technical University Eindhoven
and Universität Leipzig*

We investigate the local times of a continuous-time Markov chain on an arbitrary discrete state space. For fixed finite range of the Markov chain, we derive an explicit formula for the joint density of all local times on the range, at any fixed time. We use standard tools from the theory of stochastic processes and finite-dimensional complex calculus.

We apply this formula in the following directions: (1) we derive large deviation upper estimates for the normalized local times beyond the exponential scale, (2) we derive the upper bound in Varadhan's lemma for any measurable functional of the local times, and (3) we derive large deviation upper bounds for continuous-time simple random walk on large subboxes of $\mathbb{Z}^d$ tending to $\mathbb{Z}^d$ as time diverges. We finally discuss the relation of our density formula to the Ray–Knight theorem for continuous-time simple random walk on $\mathbb{Z}$, which is analogous to the well-known Ray–Knight description of Brownian local times.

**1. Introduction.** Let $\Lambda$ be a finite or countably infinite set and let $A = (A_{x,y})_{x,y \in \Lambda}$ be the generator, sometimes called the $Q$-matrix, of a continuous-time Markov chain $(X_t)_{t \in [0,\infty)}$ on $\Lambda$. Under the measure $\mathbb{P}_a$, the chain starts at $X_0 = a \in \Lambda$, and by $\mathbb{E}_a$ we denote the corresponding expectation. The main object of our study are the *local times*, defined by

$$(1.1) \qquad \ell_T(x) = \int_0^T \mathbb{1}_{\{X_s = x\}} \, ds, \qquad x \in \Lambda, T > 0,$$

Received October 2005; revised October 2006.
[1]Supported by the Natural Sciences and Engineering Research Council of Canada.
[2]Supported in part by Netherlands Organization for Scientific Research (NWO).
[3]Supported in part by German Science Foundation.

which register the amount of time the chain spends in $x$ up to time $T$. We have $\langle \ell_T, V \rangle = \int_0^T V(X_s)\, ds$ for any bounded function $V : \Lambda \to \mathbb{R}$, where $\langle \cdot, \cdot \rangle$ denotes the standard inner product on $\mathbb{R}^\Lambda$. We conceive the normalized local times tuple, $\frac{1}{T}\ell_T = (\frac{1}{T}\ell_T(x))_{x \in \Lambda}$, as a random element of the set $\mathcal{M}_1(\Lambda)$ of probability measures on $\Lambda$.

The local times tuple $\ell_T = (\ell_T(x))_{x \in \Lambda}$, and in particular its large-$T$ behavior, are of fundamental interest in many branches and applications of probability theory. We are particularly interested in the large deviation of $\frac{1}{T}\ell_T$. A by now classical result [12, 21] states, for a *finite* state space $\Lambda$, a large deviation principle for $\frac{1}{T}\ell_T$, for any starting point $a \in \Lambda$, on the scale $T$. More precisely, for any closed set $\Gamma \subseteq \mathcal{M}_1(\Lambda)$,

$$(1.2) \qquad \limsup_{T \to \infty} \frac{1}{T} \log \mathbb{P}_a\left(\frac{1}{T}\ell_T \in \Gamma\right) \leq -\inf_{\mu \in \Gamma} I_A(\mu),$$

and, for any open set $G \subseteq \mathcal{M}_1(\Lambda)$,

$$(1.3) \qquad \liminf_{T \to \infty} \frac{1}{T} \log \mathbb{P}_a\left(\frac{1}{T}\ell_T \in G\right) \geq -\inf_{\mu \in G} I_A(\mu).$$

The rate function $I_A$ may be written

$$(1.4) \qquad I_A(\mu) = -\inf\left\{ \left\langle Ag, \frac{\mu}{g} \right\rangle \,\Big|\, g : \Lambda \to (0, \infty) \right\}.$$

In case that $A$ is a symmetric matrix, $I_A(\mu) = \|(-A)^{1/2}\sqrt{\mu}\|_2^2$ is equal to the Dirichlet form of $A$ applied to $\sqrt{\mu}$. The topology used on $\mathcal{M}_1(\Lambda)$ is the weak topology induced by convergence of integrals against all bounded functions $\Lambda \to \mathbb{R}$, that is, the standard topology of pointwise convergence since $\Lambda$ is assumed finite. For infinite $\Lambda$, versions of this large deviations principle may be formulated for the restriction of the chain to some finite subset of $\Lambda$. A standard way of proving the above principle of large deviations is via the Gärtner–Ellis theorem; see [11] for more background on large deviation theory. One of the major corollaries is Varadhan's lemma, which states that

$$(1.5) \qquad \lim_{T \to \infty} \frac{1}{T} \log \mathbb{E}_a[e^{TF((1/T)\ell_T)}] = -\inf_{\mu \in \mathcal{M}_1(\Lambda)}[I_A(\mu) - F(\mu)],$$

for any function $F : \mathcal{M}_1(\Lambda) \to \mathbb{R}$ that is bounded and continuous in the above topology. We would like to stress that in many situations it is the upper bound in (1.5) that is difficult to prove since $F$ often fails to be upper semicontinuous. [However, often $F$ turns out to be lower semicontinuous or well approximated by lower semicontinous functions, so that the proof of the lower bound in (1.5) is often simpler.]

In the present paper, we considerably strengthen the above large deviation principle and the assertion in (1.5) by presenting an *explicit density of the random variable* $\ell_T$, that is, a joint density of the tuple $(\ell_T(x))_{x \in \Lambda}$, for any fixed $T > 0$. We do this for either a finite state space $\Lambda$ or for the restriction to a finite subset. This formula opens up several new possibilities, such as:

(1) more precise asymptotics for the probabilities in (1.2) and (1.3) and for the expectation on the left-hand side of (1.5),

(2) the validity of (1.5) for many discontinuous functions $F$,

(3) versions of the large deviation principle for rescaled versions of the local times on state spaces $\Lambda = \Lambda_T$ coupled with $T$ and growing to some infinite set.

Clearly, a closed analytical formula for the density of the local times is quite interesting in its own right. Unfortunately, our expression for the local times density is rather involved and is quite hard to evaluate asymptotically. Actually, not even the nonnegativity of the density can be easily seen from our formula. Luckily, *upper bounds* on the density are more easily obtained. We will be able to use these upper bounds to derive proofs of (1.2) and of the upper bound in (1.5) for *every measurable set* $\Gamma$, respectively, for *every measurable function* $F$, which is a great improvement.

This paper is organized as follows. In Section 2, we identify the density of the local times in Theorem 2.1, and prove Theorem 2.1. In Section 3, we use Theorem 2.1 to prove large deviation upper bounds in Theorem 3.6. Finally, we close in Section 4 by discussing our results, by relating them to the history of the problem and by discussing the relation to the Ray–Knight theorem.

**2. Density of the local times.** In this section, we present our fundamental result, Theorem 2.1, which is the basis for everything that follows. By

$$(2.1) \qquad R_T = \mathrm{supp}(\ell_T) = \{X_s : s \in [0,T]\} \subseteq \Lambda$$

we denote the *range* of the Markov chain. Note that given $\{R_T \subseteq R\}$ for some finite set $R \subseteq \Lambda$, the random tuple $(\ell_T(x))_{x \in R}$ does not have a density with respect to the Lebesgue measure, since the event $\{\ell_T(x) = 0\}$ occurs with positive probability for any $x \in R$, except for the initial site of the chain. However, given $\{R_T = R\}$ for some $R \subseteq \Lambda$, the tuple $(\ell_T(x))_{x \in R}$ takes values in the simplex

$$(2.2) \qquad \mathcal{M}_T^+(R) = \left\{ l : R \to (0, \infty) \, \Big| \, \sum_{x \in R} l(x) = T \right\},$$

which is a convex open subset of the hyperplane in $\mathbb{R}^R$ that is perpendicular to $\mathbb{1}$. It will turn out that on $\{R_T = R\}$, the tuple $(\ell_T(x))_{x \in R}$ has a density with respect to the Lebesgue measure $\sigma_T$ on $\mathcal{M}_T^+(R)$ defined by the disintegration of Lebesgue measure into surface measures,

$$(2.3) \qquad \int d^R l\, F(l) = \int_0^\infty dT \int_{\mathcal{M}_T^+} \sigma_T(dl) F(l),$$

where $F\colon (0,\infty)^R \to \mathbb{R}$ is bounded and continuous with compact support.

We need some notation. Let $R \subseteq \Lambda$ and let $a,b \in R$. For a matrix $M = (M_{x,y})_{x,y \in \Lambda}$ we denote by $\det_{ab}^{(R)}(M)$ the $(b,a)$ cofactor of the $R \times R$-submatrix of $M$, namely, the determinant of the matrix $(1_{x \neq b} M_{x,y} 1_{y \neq a} + 1_{x=b,y=a})_{x,y \in R}$. We write $\det_{ab}$ instead of $\det_{ab}^{(\Lambda)}$ when no confusion can arise. By $\partial_l$ we denote the $\Lambda \times \Lambda$-diagonal matrix with $(x,x)$-entry $\partial_{l_x}$, which is the partial derivative with respect to $l_x$. Hence, $\det_{ab}^{(R)}(M + \partial_l)$ is a linear differential operator of order $|R| - 2 + \delta_{a,b}$.

Then our main result reads as follows:

THEOREM 2.1 (Density of the local times). *Let $\Lambda$ be a finite or countably infinite set with at least two elements and let $A = (A_{x,y})_{x,y \in \Lambda}$ be the conservative generator of a continuous-time Markov chain on $\Lambda$. Fix a finite subset $R$ of $\Lambda$ and sites $a,b \in R$. Then, for every $T > 0$ and for every bounded measurable function $F\colon \mathcal{M}_T^+(R) \to \mathbb{R}$,*

$$(2.4) \qquad \mathbb{E}_a[F(\ell_T) 1_{\{X_T = b\}} 1_{\{R_T = R\}}] = \int_{\mathcal{M}_T^+(R)} F(l) \rho_{ab}^{(R)}(l) \sigma_T(dl),$$

*where, for $l \in \mathcal{M}_T^+(R)$,*

$$(2.5) \quad \rho_{ab}^{(R)}(l) = \det_{ab}^{(R)}(-A + \partial_l) \int_{[0,2\pi]^R} e^{\sum_{x,y \in R} A_{x,y} \sqrt{l_x} \sqrt{l_y} e^{\mathrm{i}(\theta_x - \theta_y)}} \prod_{x \in R} \frac{d\theta_x}{2\pi}.$$

Alternative expressions for the density $\rho_{ab}^{(R)}$ are found in Proposition 2.5 below. Note that the density $\rho_{ab}^{(R)}$ does not depend on the values of the generator outside $R$, nor on $T$. The formula for the density is explicit, but quite involved, in particular as it involves determinants of large matrices, additional multiple integrals, and various partial derivatives. For example, it is not clear from (2.5) that $\rho_{ab}^{(R)}$ is nonnegative. Nevertheless, the formula allows us to prove rather precise and transparent large deviation upper bounds for the local times as we shall see later. As we will discuss in more detail in Section 4, Theorem 2.1 finds its roots in the work of Luttinger [29] who expressed expectations of functions of the local times in terms of integrals in which there are "functions" of anticommuting differential forms (Grassman variables). It is not clear from his work that the Grassman variables can be removed without creating intractable expressions. Theorem 2.1 accomplishes this removal. We also provide a proof that makes no overt use of Grassman variables; the determinant is their legacy.

To prepare for the proof, we need the following two lemmas and some notation. We write $\phi = u + \mathrm{i}v$ and $\overline{\phi} = u - \mathrm{i}v$, where $u,v \in \mathbb{R}^\Lambda$, and we use $d^\Lambda u \, d^\Lambda v$ to denote the Lebesgue measure on $\mathbb{R}^\Lambda \times \mathbb{R}^\Lambda$. Let $\langle \phi, \psi \rangle = \sum_{x \in \Lambda} \phi_x \psi_x$ be the real inner product on $\mathbb{C}^\Lambda$.

LEMMA 2.2. *Let $\Lambda$ be a finite set, and let $M \in \mathbb{C}^{\Lambda \times \Lambda}$. If $\Re\langle \phi, M\overline{\phi}\rangle > 0$ for any $\phi \in \mathbb{C}^\Lambda \setminus \{0\}$, then*

$$\int_{\mathbb{R}^\Lambda \times \mathbb{R}^\Lambda} d^\Lambda u\, d^\Lambda v\, e^{-\langle \phi, M\overline{\phi}\rangle} = \frac{\pi^{|\Lambda|}}{\det(M)}. \tag{2.6}$$

REMARK 2.3. By introducing polar coordinates $(l, \theta) \in [0, \infty)^\Lambda \times [0, 2\pi]^\Lambda$ via

$$\phi_x = \sqrt{l_x} e^{\mathrm{i}\theta_x}, \qquad x \in \Lambda, \tag{2.7}$$

we can transform

$$d^\Lambda u\, d^\Lambda v = \pi^{|\Lambda|} \prod_{x \in \Lambda} \left( dl_x \frac{d\theta_x}{2\pi} \right) = 2^{-|\Lambda|} d^\Lambda l\, d^\Lambda \theta \tag{2.8}$$

and can rewrite (2.6) in the form

$$\int_{[0,\infty)^\Lambda \times [0,2\pi]^\Lambda} d^\Lambda l \frac{d^\Lambda \theta}{(2\pi)^{|\Lambda|}} e^{-\langle \phi, M\overline{\phi}\rangle} = \frac{1}{\det(M)}. \tag{2.9}$$

PROOF OF LEMMA 2.2. We define the complex inner product $(\phi, \psi) = \langle \phi, \overline{\psi}\rangle$. Any unitary matrix $U \in \mathbb{C}^{\Lambda \times \Lambda}$ defines a complex linear transformation on $\mathbb{C}^\Lambda$ by $\phi' = U\phi$. By writing $\phi = u + \mathrm{i}v$ and $\phi' = u' + \mathrm{i}v'$ we obtain a real linear transformation $\widetilde{U} : (u, v) \mapsto (u', v')$ on $\mathbb{R}^\Lambda \oplus \mathbb{R}^\Lambda$. The map $\widetilde{U}$ is orthogonal, because

$$\langle u', u'\rangle + \langle v', v'\rangle = (\phi', \phi') = (U\phi, U\phi) = (\phi, \phi) = \langle u, u\rangle + \langle v, v\rangle.$$

Let $M^*$ be the adjoint to $M$ so that $(\phi, M\psi) = (M^*\phi, \psi)$. First we consider the case where $M = M^*$. The hypothesis $\Re\langle \phi, M\overline{\phi}\rangle > 0$ can be rewritten as $(\phi, M\phi) > 0$, so that $M$ has throughout positive eigenvalues $\lambda_x$, $x \in \Lambda$. Since $M$ is self-adjoint there exists a unitary transformation $U$ such that $U^*MU = D$, where $D$ is diagonal with diagonal entries $D_{x,x} = \lambda_x > 0$. Thus, by the change of variables $(u', v') = \widetilde{U}(u, v)$,

$$\int d^\Lambda u\, d^\Lambda v\, e^{-(\phi, M\phi)} = \int d^\Lambda u\, d^\Lambda v\, e^{-(\phi, D\phi)}.$$

The integral on the right-hand side factors into a product of integrals

$$\prod_{x \in \Lambda} \int_\mathbb{R} du \int_\mathbb{R} dv\, e^{-\lambda_x u^2 - \lambda_x v^2} = \prod_{x \in \Lambda} \frac{\pi}{\lambda_x} = \frac{\pi^{|\Lambda|}}{\det(M)}.$$

The lemma is proved for the case $M = M^*$.

Now we turn to the case where $M^* \neq M$. Let

$$S = \frac{1}{2}(M + M^*) \quad \text{and} \quad A = \frac{1}{2\mathrm{i}}(M - M^*).$$

Thus, $S$ and $A$ are self-adjoint and $M = S + \mathrm{i}A$. Also, $(\phi, S\phi) = \Re\langle \phi, M\overline{\phi}\rangle$ which is positive by the hypothesis. Therefore the eigenvalues of $S$ are strictly positive.

For $\mu \in \mathbb{C}$ we define $M(\mu) = S + \mu A$. For $\mu$ real, the matrix $M(\mu)$ is self-adjoint. Observe that $M(\mu)$ has throughout strictly positive eigenvalues when $\mu = 0$. Hence, the real part of the characteristic polynomial of $M(\mu)$ is nonzero on $(-\infty, 0]$, and therefore bounded away from zero on $(-\infty, 0]$, for $\mu = 0$. By continuity of the real part of this polynomial in $\mu$, the latter property persists to all $\mu$ in a suitable open interval $I \subset \mathbb{R}$ containing the origin. Therefore, $M(\mu)$ has throughout strictly positive eigenvalues for all $\mu \in I$. Thus we have $(\phi, M(\mu)\phi) > 0$ for all nonzero $\phi$ and all $\mu \in I$.

Now we apply the preceding with $M = M(\mu)$, and obtain, for $\mu \in I$,

$$(2.10) \qquad \det(M(\mu)) \int d^\Lambda u \, d^\Lambda v \, e^{-(\phi, M(\mu)\phi)} = \pi^{|\Lambda|}.$$

Both sides of this equation are analytic in $\mu$ for $\Re\mu \in I$ because $\det(M(\mu))$ is a polynomial in $\mu$, and the integral of the analytic function $\exp(-(\phi, M(\mu)\phi))$ is analytic by Morera's theorem and the Fubini theorem, as well as the remark that

$$|e^{-(\phi, M(\mu)\phi)}| = |e^{-(\phi, S\phi) - \mu(\phi, A\phi)}| = e^{-(\phi, S\phi) - \Re\mu(\phi, A\phi)} = e^{-(\phi, M(\Re\mu)\phi)}.$$

By analytic continuation (2.10) holds for $\Re\mu \in I$ and in particular for $\mu = \mathrm{i}$. At $\mu = \mathrm{i}$, $M(\mu) = M$. $\square$

LEMMA 2.4. *Let $\Lambda$ be a finite set, let $M \in \mathbb{C}^{\Lambda \times \Lambda}$, and $v = (v_x)_{x \in \Lambda} \in \mathbb{C}^\Lambda$. Then, for any continuously differentiable function $g : \mathbb{C}^\Lambda \to \mathbb{R}$,*

$$(2.11) \quad \det_{ab}(M + \partial_l)(e^{\langle v, \cdot \rangle} g)(l) = e^{\langle v, l \rangle} \det_{ab}(M + V + \partial_l) g(l), \qquad l \in \mathbb{R}^\Lambda,$$

*where $V = (\delta_{xy} v_x)_{x,y \in \Lambda}$ denotes the diagonal matrix with diagonal entries $v_x$.*

PROOF. By a cofactor expansion, one sees that, for any diagonal matrix $W$, $\det_{ab}(M + W) = \sum_{X \subseteq \Lambda \setminus \{a,b\}} c_X \prod_{x \in X} W_{x,x}$ for suitable coefficients $c_X$ depending only on the entries of $M$. Analogously, $\det_{ab}(M + \partial_l) = \sum_{X \subseteq \Lambda \setminus \{a,b\}} c_X \partial_l^X$, where we used the notation $\partial_l^X = \prod_{x \in X} \partial_{l_x}$. Therefore,

$$e^{-\langle v, l \rangle} \det_{ab}(M + \partial_l)(e^{\langle v, \cdot \rangle} g)(l) = \sum_{X \subseteq \Lambda \setminus \{a,b\}} c_X e^{-\langle v, l \rangle} \partial_l^X (e^{\langle v, \cdot \rangle} g)(l)$$

$$= \sum_{X \subseteq \Lambda \setminus \{a,b\}} c_X \prod_{x \in \Lambda} (v_x + \partial_{l_x}) g(l)$$

$$= \det_{ab}(M + V + \partial_l) g(l). \qquad \square$$

PROOF OF THEOREM 2.1. We have divided the proof into six steps. In the first five steps we assume that $\Lambda$ is a finite set, and we put $R = \Lambda$. Recall the notation in Remark 2.3, which will be used throughout this proof. We abbreviate $D_l = \det_{ab}(-A - \partial_l)$.

STEP 1. *For any $v \in \mathbb{C}^\Lambda$ with $\Re v \in (-\infty, 0)^\Lambda$, for $F(l) = e^{\langle v, l \rangle}$,*

$$\text{(2.12)} \quad \int_0^\infty \mathbb{E}_a[F(\ell_T) \mathbb{1}_{\{X_T = b\}}] \, dT$$
$$= \int_{[0,\infty)^\Lambda \times [0,2\pi]^\Lambda} (D_l F)(l) e^{\langle \phi, A\overline{\phi} \rangle} \prod_{x \in \Lambda} \left( dl_x \frac{d\theta_x}{2\pi} \right).$$

PROOF. Recall that $\langle v, \ell_T \rangle = \int_0^T v(X_s) \, ds$ to obtain

$$\int_0^\infty \mathbb{E}_a[F(\ell_T) \mathbb{1}_{\{X_T = b\}}] \, dT = \int_0^\infty \mathbb{E}_a[e^{\int_0^T v(X_s) \, ds} \mathbb{1}_{\{X_T = b\}}] \, dT$$
$$\text{(2.13)} \qquad\qquad = \int_0^\infty (e^{T(A+V)})_{a,b} \, dT$$
$$= (-A - V)^{-1}_{a,b},$$

where $V$ is the diagonal matrix with $(x,x)$-entry $v_x$, and $M_{x,y}$ denotes the $(x,y)$-entry of a matrix $M$. In order to see the last identity in (2.13), we note that

$$\text{(2.14)} \quad \int_0^\infty (e^{T(A+V)})_{a,b} \, dT = \left( \int_0^\infty e^{T(A+V)} \, dT \right)_{a,b},$$

and that

$$\text{(2.15)} \quad (A+V) \int_0^\infty e^{T(A+V)} \, dT = \int_0^\infty \frac{d}{dT} e^{T(A+V)} \, dT = -I.$$

By Cramér's rule, followed by (2.9),

$$(-A - V)^{-1}_{a,b} = \frac{\det_{ab}(-A - V)}{\det(-A - V)}$$
$$\text{(2.16)} \qquad = \int \det_{ab}(-A - V) e^{\langle \phi, (A+V)\overline{\phi} \rangle} \prod_{x \in \Lambda} \left( dl_x \frac{d\theta_x}{2\pi} \right)$$
$$= \int \det_{ab}(-A - V) e^{\langle v, l \rangle} e^{\langle \phi, A\overline{\phi} \rangle} \prod_{x \in \Lambda} \left( dl_x \frac{d\theta_x}{2\pi} \right).$$

We use Lemma 2.4 with $g = 1$ and $M = A$ to obtain that

$$\det_{ab}(-A - V) e^{\langle v, l \rangle} = (-1)^{|\Lambda|-1} e^{\langle v, l \rangle} \det_{ab}(A + V)$$
$$\text{(2.17)} \qquad\qquad = (-1)^{|\Lambda|-1} \det_{ab}(A + \partial_l) e^{\langle v, l \rangle}$$
$$= \det_{ab}(-A - \partial_l) e^{\langle v, l \rangle} = (D_l F)(l),$$

where we recall that $D_l = \det_{ab}(-A - \partial_l)$. Substituting this in (2.16) and combining this with (2.13), we conclude that (2.12) holds. □

STEP 2. *The formula* (2.12) *is also valid for functions $F$ of the form*

(2.18)
$$F(l) = \prod_{x \in \Lambda} (e^{v_x l_x} f_x(l_x)), \qquad f_x \in \mathcal{C}^2((0, \infty)),$$

$$\mathrm{supp}(f_x) \subseteq (0, \infty) \text{ compact}, \qquad \Re v_x < 0.$$

PROOF. Note that (2.12) is linear in $F$ and so if we know it for exponentials, then we obtain it for linear combinations of exponentials. In more detail, consider the Fourier representation $f_x(l_x) = \int_{\mathbb{R}} \widehat{f}_x(w_x) e^{i w_x l_x} \, dw_x$. Apply (2.12) for $v$ replaced by $v + iw$ with $w \in \mathbb{R}^\Lambda$ to obtain

$$\int_0^\infty \mathbb{E}_a[e^{\langle v, \ell_T \rangle} e^{i \langle w, \ell_T \rangle} \mathbb{1}_{\{X_T = b\}}] \, dT$$
$$= \int_{[0,\infty)^\Lambda \times [0,2\pi]^\Lambda} (D_l e^{\langle v + iw, l \rangle}) e^{\langle \phi, A \overline{\phi} \rangle} \prod_{x \in \Lambda} \left( dl_x \frac{d\theta_x}{2\pi} \right).$$

Now multiply both sides with $\prod_{x \in \Lambda} \widehat{f}_x(w_x)$ and integrate over $\mathbb{R}^\Lambda$ with respect to $d^\Lambda w$. Then we apply Fubini's theorem to move the $d^\Lambda w$ integration inside. From the representation

$$\widehat{f}_x(w_x) = \frac{1}{2\pi} \int_{\mathbb{R}} f_x(l_x) e^{-i w_x l_x} \, dl_x$$

we see that $\widehat{f}_x$ is continuous by the dominated convergence theorem. Furthermore, $\widehat{f}_x$ satisfies the bound

$$|\widehat{f}_x(w_x)| = \left| \frac{1}{2\pi (i w_x)^2} \int f_x(l) \frac{d^2}{dw_x^2} e^{-i w_x l} \, dl \right|$$
$$= \frac{1}{2\pi} \frac{1}{w_x^2} \left| \int f_x''(l) e^{-i w_x l} \, dl \right| \leq \frac{1}{2\pi} \frac{1}{w_x^2} \int |f_x''(l)| \, dl.$$

Hence, all functions $w_x \mapsto \widehat{f}_x(w_x)$ are absolutely integrable, and the exponentials with $\Re v_x < 0$ make the integration over $l_x$ convergent for any $x \in R$. □

In the following we abbreviate $D_l^* = \det_{ab}(-A + \partial_l)$.

STEP 3. *For $F$ as in* (2.18),

(2.19) $$\int_0^\infty \mathbb{E}_a[F(\ell_T) \mathbb{1}_{\{X_T = b\}}] \, dT = \int F(l) D_l^* e^{\langle \phi, A \overline{\phi} \rangle} \prod_{x \in \Lambda} \left( dl_x \frac{d\theta_x}{2\pi} \right).$$

PROOF. Comparing (2.12) with this formula we see that it is enough to prove that the integration by parts formula

$$(2.20) \qquad \int (D_l F)(l)\, e^{\langle \phi, A\overline{\phi}\rangle}\, d^\Lambda l = \int F(l)(D_l^* e^{\langle \phi, A\overline{\phi}\rangle})\, d^\Lambda l$$

holds for any $\theta \in [0, 2\pi]^\Lambda$. Since $D_l = \det_{ab}(-A - \partial_l)$ is a linear differential operator which is first order in each partial derivative, it suffices to consider one integral at a time and perform the integration by parts as follows: for any $x \in \Lambda$ and any fixed $(l_y)_{y \in \Lambda \setminus \{x\}}$,

$$(2.21) \qquad \int_0^\infty (-\partial_{l_x} F(l)) e^{\langle \phi, A\overline{\phi}\rangle}\, dl_x = \int_0^\infty F(l) \partial_{l_x} e^{\langle \phi, A\overline{\phi}\rangle}\, dl_x, \qquad x \in \Lambda.$$

There are no boundary contributions because the map $l_x \mapsto F(l)$ has a compact support in $(0, \infty)$. This proves (2.19). □

STEP 4. *For any $v \in \mathbb{C}^\Lambda$,*

$$(2.22) \qquad \int_0^\infty \mathbb{E}_a[e^{\langle v, \ell_T \rangle} \mathbb{1}_{\{X_T = b\}} \mathbb{1}_{\{R_T = \Lambda\}}]\, dT = \int e^{\langle v, l \rangle} D_l^* e^{\langle \phi, A\overline{\phi}\rangle} \prod_{x \in \Lambda} \left( dl_x \frac{d\theta_x}{2\pi} \right).$$

PROOF. Let $(f_n)_{n \in \mathbb{N}}$ be a uniformly bounded sequence of smooth functions with compact support in $(0, \infty)$ such that $f_n(t) \to \mathbb{1}_{(0,\infty)}(t)$ for any $t$. Choose $F(l) = F_n(l) = \prod_{x \in \Lambda}(e^{v_x l_x} f_n(l_x))$ in (2.19) and take the limit as $n \to \infty$, interchanging the limit with the integrals using the dominated convergence theorem. Observe that $\lim_{n \to \infty} F_n(\ell_T) = e^{\langle v, \ell_T \rangle} \prod_{x \in \Lambda} \mathbb{1}_{(0,\infty)}(\ell_T(x)) = e^{\langle v, \ell_T \rangle} \mathbb{1}_{\{R_T = \Lambda\}}$ almost surely. Furthermore, $\lim_{n \to \infty} F_n(l) = e^{\langle v, l \rangle}$ almost everywhere with respect to the measure $\prod_{x \in \Lambda}(dl_x \frac{d\theta_x}{2\pi})$. Thus we obtain (2.22) in the limit of (2.19). □

STEP 5. *For all $v \in \mathbb{C}^\Lambda$,*

$$(2.23) \qquad \mathbb{E}_a[e^{\langle v, \ell_T \rangle} \mathbb{1}_{\{X_T = b\}} \mathbb{1}_{\{R_T = \Lambda\}}] = \int_{\mathcal{M}_T^+(\Lambda)} e^{\langle v, l \rangle} \rho_{ab}^{(\Lambda)}(l)\, d^\Lambda l,$$

$$T > 0, \qquad a, b \in \Lambda,$$

*where $\rho_{ab}^{(\Lambda)}(l)$ is given by (2.5).*

PROOF. Recall that $\sum_{x \in \Lambda} \ell_T(x) = T$ almost surely and that $\sum_{x \in \Lambda} l_x = T$ for $l \in \mathcal{M}_T^+(\Lambda)$. Hence, without loss of generality, we can assume that $\Re v \in (-\infty, 0)^\Lambda$, since adding a constant $C \in \mathbb{R}$ to all the $v_x$ results in adding a factor of $e^{CT}$ on both sides. In (2.22) we replace $v_x$ by $v_x - \lambda$ with $\lambda > 0$.

Then (2.22) becomes

$$\int_0^\infty e^{-\lambda T} \mathbb{E}_a[e^{\langle v, \ell_T\rangle} \mathbb{1}_{\{X_T=b\}} \mathbb{1}_{\{R_T=\Lambda\}}]\, dT$$

$$= \int e^{\langle v,l\rangle} e^{-\lambda \sum_x l_x} D_l^* e^{\langle \phi, A\overline{\phi}\rangle} \prod_{x\in\Lambda} \left(dl_x \frac{d\theta_x}{2\pi}\right)$$

(2.24)

$$= \int_{(0,\infty)^\Lambda} e^{\langle v,l\rangle} e^{-\lambda \sum_x l_x} \rho_{ab}^{(\Lambda)}(l)\, d^\Lambda l$$

$$= \int_0^\infty e^{-\lambda T} \left[\int_{\mathcal{M}_T^+(\Lambda)} e^{\langle v,l\rangle} \rho_{ab}^{(\Lambda)}(l) \sigma_T(dl)\right] dT,$$

where

(2.25) $$\rho_{ab}^{(\Lambda)}(l) = \int_{[0,2\pi]^\Lambda} D_l^* e^{\langle\phi, A\overline{\phi}\rangle} \prod_{x\in\Lambda} \frac{d\theta_x}{2\pi}, \qquad l \in (0,\infty)^\Lambda.$$

In the second equation, we have interchanged the integrations over $l$ and $\theta$ and have rewritten the $\theta$ integral using (2.5). In the third equation in (2.24), we have introduced the variable $T = \sum_x l_x$ and used (2.3).

Hence we have proved that the Laplace transforms with respect to $T$ of the two sides of (2.23) coincide. As a consequence, (2.23) holds for almost every $T > 0$. Furthermore, (2.23) even holds for all $T > 0$, since both sides are continuous. Indeed, for small $h$ we have $\langle v, \ell_{T+h}\rangle = \langle v, \ell_T\rangle$, $X_{T+h} = X_T$ and $R_{T+h} = R_T$ with high probability, which easily implies the continuity of the left-hand side of (2.23). We see that the right-hand side is continuous for $T > 0$ by using the change of variable $t = T^{-1}l$ and (2.5) to rewrite the right-hand side as an integral of a continuous function of $T, t$ on the standard simplex $\mathcal{M}_{T=1}^+(\Lambda)$. □

Now we complete the proof of the theorem.

STEP 6. *The formula (2.4) holds for any finite or countably infinite state space $\Lambda$ and any finite subset $R$ of $\Lambda$.*

PROOF. It is enough to prove (2.4) for the case $F(l) = e^{\langle v,l\rangle}$ with $\Re(v) \in (-\infty, 0)^R$ because the distribution of $(\ell_T(x))_{x\in R}$ on the event $\{R_T = R\}$ is determined by its characteristic function.

Consider the Markov chain on $R$ with conservative generator $A^{(R)} = (A^{(R)}_{x,y})_{x,y\in R}$ given by

(2.26) $$A^{(R)}_{x,y} = \begin{cases} A_{x,y}, & \text{if } x \neq y, \\ -\sum_{y\in R\setminus\{x\}} A_{x,y}, & \text{if } x = y, \end{cases}$$

and let $V^{(R)}$ be the diagonal $R \times R$ matrix with $V^{(R)}_{x,x} = \sum_{y \in \Lambda \setminus R} A_{x,y}$. Then

$$A^{(R)}_{x,y} = A_{x,y} + V^{(R)}_{x,y} \qquad \forall x, y \in R. \tag{2.27}$$

When started in $R$, the Markov chain with generator $A^{(R)}$ coincides with the original one as long as no step to a site outside $R$ is attempted. Step decisions outside $R$ are suppressed. The distribution of this chain is absolutely continuous with respect to the original one. More precisely,

$$\begin{aligned}
&\mathbb{E}_a[F(\ell_T) \mathbb{1}_{\{X_T=b\}} \mathbb{1}_{\{R_T=R\}}] \\
&\qquad = \mathbb{E}^{(R)}_a[F(\ell_T) e^{-\sum_{x \in R} \ell_T(x) V^{(R)}_{x,x}} \mathbb{1}_{\{X_T=b\}} \mathbb{1}_{\{R_T=R\}}], \\
&\qquad\qquad\qquad\qquad\qquad\qquad\qquad\qquad T > 0, a, b \in R,
\end{aligned} \tag{2.28}$$

where $\mathbb{E}^{(R)}_a$ is the expectation with respect to the Markov chain on $R$ with generator $A^{(R)}$. Applying (2.23) for this chain with $e^{\langle v, l \rangle}$ replaced by

$$F_R(l) = F(l) e^{-\sum_{x \in R} l_x V^{(R)}_{x,x}} \tag{2.29}$$

and with $\Lambda$ replaced by $R$, we obtain, writing $\partial^{(R)}_l$ for the restriction of $\partial_l$ to $R \times R$,

$$\begin{aligned}
\mathbb{E}_a[F(\ell_T) \mathbb{1}_{\{X_T=b\}} \mathbb{1}_{\{R_T=R\}}] &= \mathbb{E}^{(R)}_a[F_R(\ell_T) \mathbb{1}_{\{X_T=b\}} \mathbb{1}_{\{R_T=R\}}] \\
&= \int_{\mathcal{M}^+_T(R)} F_R(l) \rho^{(R)}_{ab}(l) \sigma_T(dl) \\
&= \int_{\mathcal{M}^+_T(R)} F(l) \widetilde{\rho}^{(\Lambda,R)}_{ab}(l) \sigma_T(dl),
\end{aligned} \tag{2.30}$$

where

$$\begin{aligned}
\widetilde{\rho}^{(\Lambda,R)}_{ab}(l) &= e^{-\sum_{x \in R} l_x V^{(R)}_{x,x}} \det_{ab}(-A^{(R)} + \partial^{(R)}_l) \\
&\quad \times \int_{[0,2\pi]^R} e^{\sum_{x,y \in R} \phi_x A^{(R)}_{x,y} \overline{\phi}_y} \prod_{x \in R} \frac{d\theta_x}{2\pi}.
\end{aligned} \tag{2.31}$$

By Lemma 2.4, followed by (2.27),

$$\begin{aligned}
\widetilde{\rho}^{(\Lambda,R)}_{ab}(l) &= \det_{ab}(-A^{(R)} - V^{(R)} + \partial^{(R)}_l) \\
&\quad \times \left[ e^{-\sum_{x \in R} l_x V^{(R)}_{x,x}} \int_{[0,2\pi]^R} e^{\sum_{x,y \in R} \phi_x A^{(R)}_{x,y} \overline{\phi}_y} \prod_{x \in R} \frac{d\theta_x}{2\pi} \right] \\
&= \det_{ab}(-A^{(R)} - V^{(R)} + \partial^{(R)}_l) \\
&\quad \times \int_{[0,2\pi]^R} e^{\sum_{x,y \in R} \phi_x A_{x,y} \overline{\phi}_y} \prod_{x \in R} \frac{d\theta_x}{2\pi}
\end{aligned} \tag{2.32}$$

$$= \det_{ab}^{(R)}(-A + \partial_l) \int_{[0,2\pi]^R} e^{\sum_{x,y \in R} \phi_x A_{x,y} \overline{\phi}_y} \prod_{x \in R} \frac{d\theta_x}{2\pi}.$$

From the definition (2.5), and using (2.7), we recognize the last line as $\rho_{ab}^{(R)}(l)$. Therefore, by combining (2.32) and (2.30) we have proved (2.4) in the theorem.  □

Now we collect some alternative expressions for the density $\rho_{ab}^{(R)}$.

PROPOSITION 2.5. *Let the assumptions of Theorem 2.1 be satisfied. Let $B = ([1 - \delta_{x,y}]A_{x,y})_{x,y \in \Lambda}$ be the off-diagonal part of $A$. Then, for any finite subset $R$ of $\Lambda$ and for any sites $a, b \in R$, and for any $l \in \mathcal{M}_T^+(R)$, the following holds:*

(i)

(2.33)
$$\rho_{ab}^{(R)}(l) = e^{\sum_{x \in R} l_x A_{x,x}} \det_{ab}^{(R)}(-B + \partial_l)$$
$$\times \int_{[0,2\pi]^R} e^{\sum_{x,y \in R} B_{x,y} \sqrt{l_x} \sqrt{l_y} e^{i(\theta_x - \theta_y)}} \prod_{x \in R} \frac{d\theta_x}{2\pi}.$$

(ii) *For any $r \in (0,\infty)^R$,*

(2.34)
$$\rho_{ab}^{(R)}(l) = e^{\sum_{x \in R} l_x A_{x,x}} \det_{ab}^{(R)}(-B + \partial_l)$$
$$\times \int_{[0,2\pi]^R} e^{\sum_{x,y \in R} r_x B_{x,y} r_y^{-1} \sqrt{l_x} \sqrt{l_y} e^{i(\theta_x - \theta_y)}} \prod_{x \in R} \frac{d\theta_x}{2\pi}.$$

(iii)

(2.35)
$$\rho_{ab}^{(R)}(l) = \int_{[0,2\pi]^R} \det_{ab}^{(R)}(-B + V_{\theta,l})$$
$$\times e^{\sum_{x,y \in R} A_{x,y} \sqrt{l_x} \sqrt{l_y} e^{i(\theta_x - \theta_y)}} \prod_{x \in R} \frac{d\theta_x}{2\pi},$$

*where $V_{\theta,l} = (\delta_{x,y} v_{\theta,l}(x))_{x \in R}$ is the diagonal matrix with entries*

(2.36) $$v_{\theta,l}(x) = \sum_{z \in R} B_{x,z} \sqrt{\frac{l_z}{l_x}} e^{i(\theta_x - \theta_z)}, \qquad x \in R.$$

The formula in (2.34) will be helpful later when we derive upper bounds on $\rho_{ab}^{(R)}(l)$ in the case that $A$ is not symmetric. The remainder of the paper does not rely on the formula in (2.35). However, we find (2.35) of independent interest, since the integral in (2.35) does not involve any derivative.

PROOF OF PROPOSITION 2.5. Formula (2.33) follows from (2.5) by using Lemma 2.4.

We now prove (2.34). Fix $r \in (0, \infty)^R$ and observe that, for any $l \in (0, \infty)^R$,

$$
\begin{aligned}
&\int_{[0,2\pi]^R} e^{\sum_{x,y \in R} B_{x,y} \sqrt{l_x} \sqrt{l_y} e^{\mathrm{i}(\theta_x - \theta_y)}} \prod_{x \in R} \frac{d\theta_x}{2\pi} \\
(2.37) \quad &= \int_{[0,2\pi]^R} e^{\sum_{x,y \in R} r_x B_{x,y} r_y^{-1} \sqrt{l_x} \sqrt{l_y} e^{\mathrm{i}(\theta_x - \theta_y)}} \prod_{x \in R} \frac{d\theta_x}{2\pi}.
\end{aligned}
$$

Indeed, substituting $e^{\mathrm{i}\theta_x} = z_x$ for $x \in R$, we can rewrite the integrals as integrals over circles in the complex plane. The integrand is analytic in $z_x \in \mathbb{C} \setminus \{0\}$. Hence, the integral is independent of the curve (as long as it is closed and winds around zero precisely once), and it is equal to the integral along the centered circle with radius $r_x$ instead of radius one. Re-substituting $r_x e^{\mathrm{i}\theta_x} = z_x$, we arrive at (2.37). Comparing to (2.33), we see that we have derived (2.34).

Finally, we prove (2.35). We use (2.34) with $r = \sqrt{l}$ and interchange $\det_{ab}^{(R)}(-B + \partial_l)$ with $\int_{[0,2\pi]^R}$ (this is justified by the analyticity of the integrand in all the $l_x$ with $x \in R$). This gives that

$$
\rho_{ab}^{(R)}(l) = e^{\sum_{x \in R} l_x A_{x,x}} \int_{[0,2\pi]^R} \det_{ab}^{(R)}(-B + \partial_l) e^{\sum_{x,y \in R} l_x B_{x,y} e^{\mathrm{i}(\theta_x - \theta_y)}} \prod_{x \in R} \frac{d\theta_x}{2\pi}.
$$

Use Lemma 2.4 with $g = 1$ to see that

$$
\det_{ab}^{(R)}(-B + \partial_l) e^{\sum_{x,y \in R} l_x B_{x,y} e^{\mathrm{i}(\theta_x - \theta_y)}} = e^{\sum_{x,y \in R} l_x B_{x,y} e^{\mathrm{i}(\theta_x - \theta_y)}} \det_{ab}^{(R)}(-B + \widetilde{V}_\theta),
$$

where $\widetilde{V}_\theta = (\delta_{x,y} \widetilde{v}_\theta(x))_{x \in R}$ is the diagonal matrix with entries $\widetilde{v}_\theta(x) = \sum_{z \in R} B_{x,z} e^{\mathrm{i}(\theta_x - \theta_z)}$.

Now we use the same transformation as in (2.37): We interpret the integrals over $\theta_x$ as integrals over circles of radius $\sqrt{l_x}$ and replace them by integrals over circles with radius one. By this transformation, $\widetilde{V}_\theta$ is transformed into $V_{\theta,l}$, and the term $e^{\sum_{x,y \in R} l_x B_{x,y} e^{\mathrm{i}(\theta_x - \theta_y)}}$ is transformed into $e^{\sum_{x,y \in R} \sqrt{l_x} B_{x,y} \sqrt{l_y} e^{\mathrm{i}(\theta_x - \theta_y)}}$. Recalling that $B$ is the off-diagonal part of $A$, (2.35) follows. □

**3. Large deviation upper bounds for the local times.** In this section we use Theorem 2.1 to derive sharp upper bounds for the probability in (1.2) and for the expectation in (1.5) for fixed $T$ and fixed finite ranges of the local times. The main term in this estimate is given in terms of the rate function $I_A$. The main value of our formula, however, comes from the facts that (1) the error term is controlled on a subexponential scale, (2) the set $\Gamma$ in (1.2) is just assumed measurable, and (3) the functional $F$ in (1.5) is just assumed measurable. Let us stress that this formula is extremely useful, since the functional $F$ is not upper semicontinuous nor bounded in many important applications.

In Section 3.1 we give a pointwise upper bound for the density, in Section 3.2 we apply it to derive upper bounds for the probability in (1.2) and for the expectation in (1.5), and in Section 3.3 we consider the same problem for state spaces $\Lambda = \Lambda_T \subseteq \mathbb{Z}^d$ depending on $T$ and increasing to $\mathbb{Z}^d$.

3.1. *Pointwise upper bound for the density.* Here is a pointwise upper bound for the density. Recall the rate function $I_A$ introduced in (1.4).

PROPOSITION 3.1 (Upper bound for $\rho_{ab}^{(R)}$). *Under the assumptions of Theorem* 2.1, *for any finite subset $R$ of $\Lambda$, and for any $a, b \in R$, any $T > 0$ and any $l \in \mathcal{M}_T^+(R)$,*

$$
\rho_{ab}^{(R)}(l) \leq e^{-TI_A((1/T)l)} \Bigg( \prod_{x \in R \setminus \{a,b\}} \sqrt{\frac{T}{l_x}} \Bigg) \eta_R^{|R|-1}
$$
(3.1)
$$
\times e^{[\eta_R^{-1} + (4\eta_R^2 T)^{-1}] \sum_{x,y \in R} \sqrt{l_x} g_y B_{x,y}/(\sqrt{l_y} g_x)},
$$

*where $g \in (0, \infty)^R$ is a minimizer in* (1.4) *for $\mu = \frac{l}{T}$ and*

(3.2) $$\eta_R = \max\Bigg\{ \max_{x \in R} \sum_{y \in R \setminus \{x\}} |B_{x,y}|, \max_{y \in R} \sum_{x \in R \setminus \{y\}} |B_{x,y}|, 1 \Bigg\},$$

*where $B = ([1 - \delta_{x,y}]A_{x,y})_{x,y \in \Lambda}$ is the off-diagonal part of $A$.*

REMARK 3.2. If $A$ (and hence $B$) is symmetric, then $g = \sqrt{\mu}$ is the minimizer in (1.4), and we have $I_A(\mu) = \|(-A)^{1/2}\sqrt{\mu}\|_2^2$. In this case the upper bound simplifies to

(3.3) $$\rho_{ab}^{(R)}(l) \leq e^{-TI_A((1/T)l)} \Bigg( \prod_{x \in R \setminus \{a,b\}} \sqrt{\frac{T}{l_x}} \Bigg) \eta_R^{|R|-1} e^{|R|[1+(4\eta_R T)^{-1}]}.$$

The proof of Proposition 3.1 makes use of three lemmas that we will state and prove first.

LEMMA 3.3. *Let $\widetilde{B} \in [0, \infty)^{R \times R}$ be any matrix with nonnegative elements, and let $Q \subseteq R$. Then*

$$
0 \leq \partial_l^Q \int_{[0,2\pi]^R} e^{\sum_{x,y \in R} \widetilde{B}_{x,y} \sqrt{l_x} \sqrt{l_y} e^{i(\theta_x - \theta_y)}} \prod_{x \in R} \frac{d\theta_x}{2\pi} \leq \partial_l^Q e^{\sum_{x,y \in R} \widetilde{B}_{x,y} \sqrt{l_x} \sqrt{l_y}},
$$
(3.4)
$$
l \in (0, \infty)^R,
$$

where $\partial_l^Q = \prod_{x \in Q} \partial_{l_x}$.

PROOF. Write $e^{\sum_{x,y \in R} \cdots} = \prod_{x,y \in R} e^{\cdots}$ and expand the exponentials as power series. For $n = (n_{x,y})_{x,y \in R} \in \mathbb{N}_0^{R \times R}$, we write $n! = \prod_{x,y \in R} n_{x,y}!$. Then we obtain

$$\partial_l^Q \int_{[0,2\pi]^R} e^{\sum_{x,y \in R} \widetilde{B}_{x,y} \sqrt{l_x} \sqrt{l_y} e^{i(\theta_x - \theta_y)}} \prod_{x \in R} \frac{d\theta_x}{2\pi}$$

$$(3.5) \quad = \sum_{n \in \mathbb{N}_0^{R \times R}} \frac{1}{n!} \partial_l^Q \left[ \prod_{x,y \in R} (\widetilde{B}_{x,y} \sqrt{l_x} \sqrt{l_y})^{n_{x,y}} \right.$$

$$\left. \times \int_{[0,2\pi]^R} e^{i \sum_{x,y \in R} n_{x,y}(\theta_x - \theta_y)} \prod_{x \in R} \frac{d\theta_x}{2\pi} \right].$$

After rewriting the exponent in the integral on the right-hand side using $\sum_{x,y} n_{x,y}(\theta_x - \theta_y) = \sum_x n_x \theta_x$, where $n_x = \sum_y (n_{x,y} - n_{y,x})$, it is clear that the integral equals one or zero. Hence, the lower bound in (3.4) is clear, and the upper bound comes from replacing the integral by one and a resummation over $n$. $\square$

LEMMA 3.4. *Fix any matrix $B \in \mathbb{R}^{R \times R}$, let $a, b \in R$, and let $f : (0, \infty)^R \to \mathbb{R}$ be any function with nonnegative derivatives, that is, $\partial_l^Q f(l) \geq 0$ for all $Q \subseteq R$. Then*

$$(3.6) \quad |\det_{ab}^{(R)}(-B + \partial_l) f| \leq \eta_R \prod_{x \in R \setminus \{a,b\}} (\eta_R + \partial_{l_x}) f,$$

*where $\eta_R$ is defined in (3.2).*

PROOF. Recalling that the determinant is the (signed) volume subtended by the rows, we can bound a determinant by the product of the lengths of the rows. This is called the Hadamard bound and it applies to any real square matrix. Therefore, for $X \subseteq R$ and $a, b \in X$,

$$|\det_{ab}^{(X)}(-B)| \leq \prod_{x \in X \setminus \{b\}} \|B_x\| \leq \prod_{x \in X \setminus \{b\}} \eta_R = \eta_R^{|X|-1},$$

where $B_x$ is the row $x$ of $B$ after eliminating the $a$th column, and $\|\cdot\|$ is the Euclidean length, which is bounded by $\eta_R$ because $\sum |a_i|^2 \leq (\sum |a_i|)^2$. Also,

$$(3.7) \quad \det_{ab}^{(R)}(-B + \partial_l) f(l)$$
$$= \sum_{\sigma : R \setminus \{b\} \to R \setminus \{a\}} \text{sign}(\hat{\sigma}) \prod_{x \in R \setminus \{a\}} (-B_{x,\sigma_x} + \delta_{x,\sigma_x} \partial_{l_x}) f(l),$$

where the sum over $\sigma$ is over all bijections $R \setminus \{b\} \to R \setminus \{a\}$, and where $\operatorname{sign}(\hat\sigma)$ is the sign of the permutation $\hat\sigma \colon R \mapsto R$ obtained by letting $\hat\sigma_x = \sigma_x$ for $x \neq b$ and $\hat\sigma_b = a$. Expanding the product, we obtain

$$
\begin{aligned}
&\det_{ab}^{(R)}(-B+\partial_l)f(l) \\
&= \sum_{Q \subseteq R\setminus\{a,b\}} \sum_{\sigma: Q^c\setminus\{b\} \to Q^c\setminus\{a\}} \operatorname{sign}(\hat\sigma)\left(\prod_{x \in Q^c\setminus\{b\}} (-B_{x,\sigma_x})\right) \\
&\qquad\qquad\qquad\qquad\qquad\qquad \times \left(\prod_{x \in Q} \partial_{l_x}\right) f(l) \\
&= \sum_{Q \subseteq R\setminus\{a,b\}} \det_{ab}^{(Q^c)}(-B)\left(\prod_{x \in Q} \partial_{l_x}\right) f(l),
\end{aligned}
\tag{3.8}
$$

where we write $Q^c = R \setminus Q$. Take absolute values and bound the cofactor using the Hadamard bound,

$$
\begin{aligned}
&|\det_{ab}^{(R)}(-B+\partial_l)f(l)| \\
&\leq \sum_{Q \subseteq R\setminus\{a,b\}} \eta_R^{|Q^c\setminus\{b\}|}\left(\prod_{x \in Q} \partial_{l_x}\right) f(l) \\
&= \eta_R \sum_{Q \subseteq R\setminus\{a,b\}} \left(\prod_{x \in (R\setminus\{a,b\})\setminus Q} \eta_R\right)\left(\prod_{x \in Q} \partial_{l_x}\right) f(l) \\
&= \eta_R \prod_{x \in R\setminus\{a,b\}} (\eta_R + \partial_{l_x}) f(l). \qquad\square
\end{aligned}
\tag{3.9}
$$

LEMMA 3.5. *Fix any finite subset $R$ of $\Lambda$, let $\widetilde{B} \in [0,\infty)^{R \times R}$ be any matrix with nonnegative elements, and fix $a,b \in R$. Then, for any $T > 0$ and any $l \in \mathcal{M}_T^+$,*

$$
\det_{ab}^{(R)}(-B+\partial_l)\int_{[0,2\pi]^R} e^{\sum_{x,y\in R}\widetilde{B}_{x,y}\sqrt{l_x}\sqrt{l_y}e^{\mathrm{i}(\theta_x-\theta_y)}} \prod_{x \in R} \frac{d\theta_x}{2\pi}
$$

$$
\leq e^{\sum_{x,y\in R}\widetilde{B}_{x,x}\sqrt{l_x}\sqrt{l_y}}\left(\prod_{x\in R\setminus\{a,b\}}\sqrt{\frac{T}{l_x}}\right)\eta_R^{|R|-1}
$$

$$
\times e^{[\eta_R^{-1}+(4\eta_R^2 T)^{-1}]\sum_{x,y\in R}\widetilde{B}_{x,y}},
\tag{3.10}
$$

*where $\eta_R$ is defined in* (3.2).

PROOF. By Lemma 3.4 followed by Lemma 3.3, we obtain

$$\text{(3.11)} \quad \text{l.h.s. of (3.10)} \leq \eta_R \prod_{x \in R \setminus \{a,b\}} (\eta_R + \partial_{l_x}) e^{\sum_{x,y \in R} \widetilde{B}_{x,y} \sqrt{l_x} \sqrt{l_y}}.$$

Substitute $t_x = \frac{\sqrt{l_x}}{\sqrt{T}} \in [0,1]$ and abbreviate $f(t) = e^{T \sum_{x,y \in R} \widetilde{B}_{x,y} t_x t_y}$. By the chain rule, $\partial_{l_x} = \frac{1}{2T} \frac{1}{t_x} \partial_{t_x}$. Then

$$\text{(3.12)} \quad \begin{aligned} \text{l.h.s. of (3.10)} &\leq \eta_R^{|R|-1} \prod_{x \in R \setminus \{a,b\}} \left(1 + \frac{1}{2\eta_R T} \frac{1}{t_x} \partial_{t_x}\right) f(t) \\ &\leq \eta_R^{|R|-1} \left(\prod_{x \in R \setminus \{a,b\}} \frac{1}{t_x}\right) \prod_{x \in R \setminus \{a,b\}} \left(1 + \frac{1}{2\eta_R T} \partial_{t_x}\right) f(t), \end{aligned}$$

where we have used that $t_x \leq 1$. Since *all* $t$ derivatives (not just the first order derivatives) of $f$ are nonnegative since $\widetilde{B}_{x,y} \geq 0$, we can add in some extra derivatives and continue the bound with

$$\text{(3.13)} \quad \begin{aligned} \text{l.h.s. of (3.10)} &\leq \eta_R^{|R|-1} \left(\prod_{x \in R \setminus \{a,b\}} \frac{1}{t_x}\right) \\ &\quad \times \prod_{x \in R \setminus \{a,b\}} \left(\sum_{n=0}^{\infty} \frac{1}{n!} \frac{\partial_{t_x}^n}{(2\eta_R T)^n}\right) f(t) \\ &= \eta_R^{|R|-1} \left(\prod_{x \in R \setminus \{a,b\}} \frac{1}{t_x}\right) f(t + (2\eta_R T)^{-1} \mathbb{1}_R), \end{aligned}$$

where the last equation follows from Taylor's theorem, and $\mathbb{1}_R : R \to \{1\}$ is the constant function.

Recalling that $t_x \leq 1$, we may estimate

$$\frac{1}{T} \log f(t + (2\eta_R T)^{-1} \mathbb{1}_R)$$

$$= \sum_{x,y \in R} \widetilde{B}_{x,y} t_x t_y + \frac{1}{2\eta_R T} \sum_{x,y \in R} \widetilde{B}_{x,y}(t_x + t_y) + \frac{1}{(2\eta_R T)^2} \sum_{x,y \in R} \widetilde{B}_{x,y}$$

$$\leq \frac{1}{T} \log f(t) + \frac{1}{T} \left[\frac{1}{\eta_R} + \frac{1}{4\eta_R^2 T}\right] \sum_{x,y \in R} \widetilde{B}_{x,y}.$$

We conclude that

$$\text{(3.14) l.h.s. of (3.10)} \leq \eta_R^{|R|-1} \left(\prod_{x \in R \setminus \{a,b\}} \frac{1}{t_x}\right) f(t) e^{[\eta_R^{-1} + (4\eta_R^2 T)^{-1}] \sum_{x,y \in R} \widetilde{B}_{x,y}}.$$

Resubstituting $t_x = \sqrt{l_x/T}$ and $f(t) = e^{T\sum_{x,y\in R}\widetilde{B}_{x,y}t_xt_y}$, the lemma is proved. □

PROOF OF PROPOSITION 3.1. Fix any $r \in (0,\infty)$ and recall the representation of the density $\rho_{ab}^{(R)}$ in (2.34). Now apply Lemma 3.5 for $\widetilde{B} = (r_x B_{x,y} r_y^{-1})_{x,y\in R}$, to obtain

$$\rho_{ab}^{(R)}(l) \leq e^{\sum_{x,y\in R} r_x\sqrt{l_x}A_{x,y}\sqrt{l_y}r_y^{-1}} \left(\prod_{x\in R\setminus\{a,b\}} \sqrt{\frac{T}{l_x}}\right) \eta_R^{|R|-1} e^{[\eta_R^{-1}+(4\eta_R^2 T)^{-1}]\sum_{x,y\in R}\widetilde{B}_{x,y}}.$$

Now we choose $r = \sqrt{l}/g$, where $g \in (0,\infty)^R$ is a minimizer in (1.4) for $\mu = \frac{1}{T}l$. This implies the bound in (3.1). □

3.2. *Upper bounds in the LDP and in Varadhan's lemma.* In this section we specialize to Markov chains having a symmetric generator $A$ and give a simple upper bound for the left hand side of (1.2) and for the expectation in (1.5). Recall from the text below (1.4) that, in the present case of a symmetric generator, $I_A(\mu) = \|(-A)^{1/2}\sqrt{\mu}\|_2^2$ for any probability measure $\mu$ on $\Lambda$.

THEOREM 3.6 (Large deviation upper bounds for the local times). *Let the assumptions of Theorem 2.1 be satisfied. Assume that $A$ is symmetric. Fix a finite subset $S$ of $\Lambda$. Then, for any $T \geq 1$ and any $a \in S$, with $\eta_S$ as in (3.2), the following bounds hold:*

(i) *For every measurable $\Gamma \subseteq \mathcal{M}_1(S)$,*

$$\log \mathbb{P}_a\left(\frac{1}{T}\ell_T \in \Gamma, R_T \subseteq S\right)$$
(3.15)
$$\leq -T\inf_{\mu\in\Gamma}\|(-A)^{1/2}\sqrt{\mu}\|_2^2 + |S|\log(\eta_S\sqrt{8e}T) + \log|S| + \frac{|S|}{4T}.$$

(ii) *For every measurable functional $F\colon \mathcal{M}_1(S) \to \mathbb{R}$,*

$$\log \mathbb{E}_a[e^{TF((1/T)\ell_T)}\mathbb{1}_{\{R_T\subseteq S\}}] \leq T\sup_{\mu\in\mathcal{M}_1(S)}[F(\mu) - \|(-A)^{1/2}\sqrt{\mu}\|_2^2]$$
(3.16)
$$+ |S|\log(\eta_S\sqrt{8e}T) + \log|S| + \frac{|S|}{4T}.$$

Theorem 3.6 is a significant improvement over the standard estimates known in large deviation theory. In fact, one standard technique to derive upper bounds for the left-hand side of (3.15) is the use of the exponential Chebyshev inequality and a compactness argument if $\Gamma$ is assumed closed.

One important ingredient there is a good control on the logarithmic asymptotics of the expectation in (3.16) for linear functions $F$. This technique produces an error of order $e^{o(T)}$, which can in general not be controlled on a smaller scale.

The standard technique to derive improved bounds on the expectation in (3.16) for fixed $T$ is restricted to *linear* functions $F$, say $F(\cdot) = \langle V, \cdot \rangle$. This technique goes via an eigenvalue expansion for the operator $A+V$ in the set $S$ with zero boundary condition. The main steps are the use of the Rayleigh–Ritz principle for the identification of the principal eigenvalue, and Parseval's identity. This gives basically the same result as in (3.16), but is strictly limited to linear functions $F$.

PROOF OF THEOREM 3.6. It is clear that (ii) follows from (i), hence we only prove (i).

According to Theorem 2.1, we may express the probability on the left-hand side of (3.15) as

$$(3.17) \quad \mathbb{P}_a\left(\frac{1}{T}\ell_T \in \Gamma, R_T \subseteq S\right) = \sum_{b \in S} \sum_{R \subseteq S \,:\, a,b \in R} \int_{\mathcal{M}_T^+(R) \cap \Gamma_{T,R}} \rho_{ab}^{(R)}(l)\, \sigma_T(dl),$$

where $\Gamma_{T,R} = T\Gamma_R$, and $\Gamma_R$ is the set of the restrictions of all the elements of $\Gamma$ to $R$.

We fix $a, b \in S$ and $R \subseteq S$ with $a, b \in R$ and use the bound in Proposition 3.1, more precisely, the one in (3.3). Hence, for $l \in \mathcal{M}_T^+(R) \cap \Gamma_{T,R}$, we obtain, after a substitution $l = T\mu$ in the exponent, that

$$\rho_{ab}^{(R)}(l) \leq e^{-T \inf_{\mu \in \Gamma \,:\, \operatorname{supp}(\mu) \subseteq R} \|(-A)^{1/2}\sqrt{\mu}\|_2^2} \left(\prod_{x \in R \setminus \{a,b\}} \sqrt{\frac{T}{l_x}}\right) \eta_R^{|R|-1}$$

(3.18)
$$\times e^{|R|[1+(4\eta_R T)^{-1}]}.$$

Substituting this in (3.17) and integrating over $l \in \mathcal{M}_T^+(R)$, we obtain

$$\mathbb{P}_a\left(\frac{1}{T}\ell_T \in \Gamma, R_T \subseteq S\right)$$

$$(3.19) \qquad \leq e^{-T \inf_{\mu \in \Gamma} \|(-A)^{1/2}\sqrt{\mu}\|_2^2} \eta_R^{|R|-1} e^{|R|[1+(4\eta_R T)^{-1}]}$$

$$\times \sum_{b \in S} \sum_{R \subseteq S \,:\, a,b \in R} \int_{\mathcal{M}_T^+(R)} \prod_{x \in R \setminus \{a\}} \sqrt{\frac{T}{l_x}}\, \sigma_T(dl)$$

$$\leq e^{-T \inf_{\mu \in \Gamma} \|(-A)^{1/2}\sqrt{\mu}\|_2^2} \eta_R^{|R|-1} e^{|S|[1+(4\eta_R T)^{-1}]} |S| \sqrt{8}^{|S|} T^{|S|-1}.$$

In the last integral, we have eliminated $l_a = T - \sum_{y \in R \setminus \{a\}} l_y$, have extended the $(|R|-1)$ single integration areas to $(0, T)$ and used that $\int_0^T l_x^{-1/2}\, dl_x =$ $\sqrt{2T}$. Now we use that $\eta_R$ is increasing in $R$ and greater than or equal to one to arrive at (3.15). This completes the proof of (i).  □

3.3. *Rescaled local times.*  As an application of Theorem 3.6, we now consider continuous-time simple random walk restricted to a large $T$-dependent subset $\Lambda = \Lambda_T$ of $\mathbb{Z}^d$ increasing to $\mathbb{Z}^d$. We derive the sharp upper bound in the large deviation principle for its *rescaled* local times. Assume, for some scale function $T \mapsto \alpha_T \in (0, \infty)$, that $\Lambda_T$ is equal to the box $[-R\alpha_T, R\alpha_T]^d \cap \mathbb{Z}^d$, where the scale function $\alpha_T$ satisfies

(3.20) $$1 \ll \alpha_T \ll \left(\frac{T}{\log T}\right)^{1/(d+2)} \qquad \text{as } T \to \infty.$$

We introduce the rescaled version of the local times,

$$L_T(x) = \frac{\alpha_T^d}{T} \ell_T(\lfloor \alpha_T x \rfloor), \qquad x \in \mathbb{R}^d.$$

Note that $L_T$ is a random step function on $\mathbb{R}^d$. In fact, it is a random probability density on $\mathbb{R}^d$. Its support is contained in the cube $[-R, R]^d$ if and only if the support of $\ell_T$ is contained in the box $[-R\alpha_T, R\alpha_T]^d \cap \mathbb{Z}^d$.

It is known that, as $T \to \infty$, the family $(L_T)_{T>0}$ satisfies a large deviation principle under the subprobability measures $\mathcal{P}(\cdot \cap \{\mathrm{supp}(L_T) \subseteq [-R, R]^d\})$ for any $R > 0$. The speed is $T\alpha_T^{-2}$, and the rate function is the energy functional, that is, the map $g^2 \mapsto \frac{1}{2}\|\nabla g\|_2^2$, restricted to the set of squares $g^2$ of $L^2$-normalized functions $g$ such that $g$ lies in $H^1(\mathbb{R}^d)$ and has its support in $[-R, R]^d$. The topology is the one which is induced by all the test integrals of $g^2$ against continuous and bounded functions. This large-deviation principle is proved in [20] for the discrete-time random walk, and the proof for continuous-time walks is rather similar (see also [24], where the proof of this fact is sketched). Hence, Varadhan's lemma yields precise logarithmic asymptotics for all exponential functionals of $L_T$ that are bounded and continuous in the above mentioned topology.

Note that this large deviations principle for $L_T$ is almost the same as the one which is satisfied by the normalized Brownian occupation times measures (see [12, 21]), the main difference being the speed (which is $T$ in [12, 21] instead of $T\alpha_T^{-2}$ here) and the fact that $L_T$ does not take values in the set of continuous functions $\mathbb{R}^d \to [0, \infty)$.

Here we want to point out that Theorem 3.6 yields a new method to derive upper bounds for many exponential functionals of $L_T$. For a cube $Q \subset \mathbb{R}^d$, we denote by $M_1(Q)$ the set of all probability densities $Q \to [0, \infty)$.

THEOREM 3.7.  *Fix $R > 0$, denote $Q_R = [-R, R]^d$ and fix a measurable function $F : M_1(Q_R) \to \mathbb{R}$. Introduce*

(3.21)  $\chi = \inf\{\frac{1}{2}\|\nabla g\|_2^2 - F(g^2) : g \in H^1(\mathbb{R}^d), \|g\|_2 = 1, \mathrm{supp}(g) \subseteq Q_R\}.$

*Then*

$$(3.22) \quad \limsup_{T\to\infty} \frac{\alpha_T^2}{T} \log \mathbb{E}_0\left[\exp\left\{\frac{T}{\alpha_T^2} F(L_T)\right\} \mathbb{1}_{\{\mathrm{supp}(L_T)\subseteq Q_R\}}\right] \leq -\chi,$$

*provided that*

$$(3.23) \quad \liminf_{T\uparrow\infty} \inf_{\mu\in\mathcal{M}_1(B_{R\alpha_T})} \left(\alpha_T^2 \tfrac{1}{2} \sum_{x\sim y}(\sqrt{\mu(x)}-\sqrt{\mu(y)})^2 - F(\alpha_T^d\mu(\lfloor\cdot\alpha_T\rfloor))\right)$$
$$\geq \chi.$$

PROOF. Introduce

$$F_T(\mu) = \frac{1}{\alpha_T^2} F(\alpha_T^d \mu(\lfloor\cdot\alpha_T\rfloor)), \qquad \mu \in \mathcal{M}_1(\mathbb{Z}^d),$$

then we have $\frac{1}{\alpha_T^2} F(L_T) = F_T(\frac{1}{T}\ell_T)$. Hence, Theorem 3.6(ii) yields that

$$\mathbb{E}_0\left[\exp\left\{\frac{T}{\alpha_T^2} F(L_T)\right\} \mathbb{1}_{\{\mathrm{supp}(L_T)\subseteq Q_R\}}\right] = \mathbb{E}_0[\exp\{TF_T(1/T\ell_T)\}\mathbb{1}_{\{\mathrm{supp}(\ell_T)\subseteq Q_{R\alpha_T}\}}]$$
$$\leq e^{o(T\alpha_T^{-2})} e^{-T\chi_T},$$

where

$$\chi_T = \inf_{\mu\in\mathcal{M}_1(Q_{R\alpha_T}\cap\mathbb{Z}^d)} \left(\tfrac{1}{2}\sum_{x\sim y}(\sqrt{\mu(x)}-\sqrt{\mu(y)})^2 - F_T(\mu)\right).$$

Here we used that the two error terms on the right-hand side of (3.16) are $e^{o(T\alpha_T^{-2})}$ since $\eta_S \leq 2d$ for any $S \subseteq \mathbb{Z}^d$ and because of our growth assumption in (3.20). Now (3.22) follows from (3.23). □

Theorem 3.7 proved extremely useful in the study of the parabolic Anderson model in [24]. Indeed, it was crucial in that paper to find the precise upper bound of the left hand side of (3.22) for the functional

$$F(g^2) = \int_{Q_R} g^2(x) \log g^2(x)\,dx,$$

which has bad continuity properties in the topology in which the above mentioned large deviations principle holds. However, Theorem 3.7 turned out to be applicable since the crucial prerequisite in (3.23) had been earlier provided in [22]. The main methods there were equicontinuity, uniform integrability and Arzela–Ascoli's theorem.

In the same paper [24], the functional

$$F(g^2) = -\int_{Q_R} |g(x)|^{2\gamma}\,dx \qquad \text{with some } \gamma\in(0,1),$$

was also considered. This problem arose in the study of the parabolic Anderson model for another type of potential distribution which was earlier studied in [4]. The prerequisite in (3.23) was provided in [24] using techniques from Gamma-convergence; see [2] for these techniques.

**4. Discussion.** In this section, we give some comments on the history of the problem addressed in the present paper.

4.1. *Historical background.* The formulas in this paper have been motivated by the work of the theoretical physicist J. M. Luttinger [29] who gave a (nonrigorous) asymptotic evaluation of certain path integrals. Luttinger claimed that there is an asymptotic series

$$\mathbb{E}_0[e^{-TF(\ell_T/T)}] \sim \sqrt{T}e^{-c_0T}\left(c_1 + \frac{c_2}{T} + \frac{c_3}{T^2} + \cdots\right)$$

for Brownian local times. He provided an algorithm to compute all the coefficients. He showed that his algorithm gives the Donsker–Varadhan large deviations formula for $c_0$ and he explicitly computed the central limit correction $c_1$.

In [10] Brydges and Muñoz-Maya used Luttinger's methods to verify that his asymptotic expansion is valid to all orders for a Markov process with symmetric generator and finite state space. The hypotheses are that $F$ is smooth and the variational principle that gives the large deviations coefficient $c_0$ is nondegenerate. Luttinger implicitly relies on similar assumptions when he uses the Feynman expansion for his functional integral.

Thus there remains the open problem to prove that Luttinger's series is asymptotic for more general state spaces, in particular, for Brownian motion. As far as we know, the best progress to date is in [5] where compact state spaces were considered and the asymptotics including the $c_1$ correction was verified.

Luttinger's paper used a calculus called *Grassman integration*. The background to this is that the Feynman–Kac formula provides a probabilistic representation for the propagation of elementary particles that satisfy "Bose statistics." To obtain a similar representation for elementary particles that satisfy "Fermi statistics" one is led in [3] to an analogue of integration defined as a linear functional on a non-Abelian Grassman algebra in place of the Abelian algebra of measurable functions: this is Grassman integration. An important part of this line of thought concerns a case where there is a relation called *supersymmetry*. This background gives no hint that Grassman integrals are relevant for ordinary Markov processes, but, nevertheless, Parisi and Sourlas [33] and McKane [32] noted that random walk expectations can be expressed in terms of the Grassman extension of Gaussian integration. Luttinger followed up on these papers by being much more explicit and precise about the supersymmetric representation in terms of Grassman integration and by deriving his series.

In [28] Le Jan pointed out that Grassman integration in this context is actually just ordinary integration in the context of differential forms. The differential forms are the non-Abelian algebra and the standard definition of integration of differential forms provides the linear functional. Since integration over differential forms is defined in terms of ordinary integration one can remove the differential forms, as we have done in this paper, but this obscures the underlying mechanism of supersymmetry. The formalism with differential forms is explained in [9], page 551, where it is used to study Green's function of a self-repelling walk on a hierarchical lattice. Two other applications of the same formalism are the proof of the Matrix–Tree theorem in [1] and a result on self-avoiding trees given in [8].

Luttinger found an instance of a relation between the local time of a Markov process on a state space $E$ and the square of a *Gaussian field* indexed by $E$. The first appearance of such a relation was given by Symanzik in [37]. His statement is that the sum of the local times of an ensemble of Brownian loops is the square of a Gaussian field. The references given above to Parisi–Sourlas, McKane and Luttinger removed the need for an ensemble by bringing, in its place, Grassman integration. The paper of Symanzik was not immediately rigorous because he claimed his result for Brownian motion but it makes almost immediate sense for Markov processes on finite state spaces only. Based on this work a rigorous relation between the square of a Gaussian field and local time of a random walk on a lattice was given by Brydges, Fröhlich and Spencer in [6]. Dynkin [13, 14, 15] showed that the identities of that paper can be extended to Brownian motion in one and two dimensions. In this form, the *Dynkin Isomorphism*, it became a useful tool for studying local time of diffusions and much work has been done by Rosen and Marcus in exploiting and extending these ideas, for example, see [19, 31]. The relation between the local time and the square of a Gaussian field is concealed in this paper in (2.12) which relates the local time $\ell$ to $l = |\phi|^2$ where $\phi$ is Gaussian. This is more obvious when $\phi$ is expressed as $\phi = u + iv$ instead of in terms of polar coordinates $\phi = \sqrt{l} e^{i\theta}$.

4.2. *Relation to the Ray–Knight theorem.* Our density formula in Theorem 2.1 can also be used to prove a version of the *Ray–Knight theorem* for continuous-time simple random walk on $\mathbb{Z}$. The well-known Ray–Knight theorem for one-dimensional Brownian motion (see [35], Sections XI.1-2, [25], Sections 6.3-4) was originally proved in [27, 34]. It describes the Brownian local times, observed at certain stopping times, as a homogeneous Markov chain in the spatial parameter. Numerous deeper investigations of this idea have been made, for example, for general symmetric Markov processes [19], for diffusions with fixed birth and death points on planar cycle-free graphs [17, 18], and on the relations to Dynkin's isomorphism [16, 36].

The (time and space) discrete version of the Ray–Knight theorem, that is, for simple random walk on $\mathbb{Z}$, was also introduced in [27], however it turned out there that it is not the local times on the sites, but on the edges that enjoys a Markov property. This idea has been used or reinvented a couple of times, for example, for applications to random walk in random environment [26], to reinforced random walk [38], and to random polymer measures [23].

In the present situation of continuous time and discrete space, it turns out that the local times themselves form a nice Markov chain. However, a proof appears to be missing. In fact, up to our best knowledge, [30] is the only paper that provides (the outline of) a proof, but only for the special case where the walk starts and ends in the same point.

We state the result here, but omit the proof. The proof will appear in an extended version [7]. We first introduce some notation. For fixed $b \in \mathbb{Z}$, we denote

$$(4.1) \qquad T_b^h = \inf\{t > 0 : \ell_t(b) > h\}, \qquad h > 0,$$

the right-continuous inverse of the map $t \mapsto \ell_t(b)$. We denote by

$$(4.2) \qquad I_0(h) = \sum_{i=0}^{\infty} \frac{h^{2i}}{2^i (i!)^2},$$

the modified Bessel function.

THEOREM 4.1 (Ray–Knight theorem for continuous-time random walks). *Let $\ell_T$ defined in (1.1) be the local times of continuous-time simple random walk $(X_t)_{t>0}$ on $\mathbb{Z}$. Let $b \in \mathbb{N}$ and $h > 0$.*

(i) *Under $\mathbb{P}_0$, the process $(\ell_{T_b^h}(b-x))_{x=0}^b$ is a time-homogeneous discrete-time Markov chain on $(0, \infty)$, starting at $h$, with transition density given by*

$$(4.3) \qquad f(h_1, h_2) = e^{-h_1 - h_2} I_0(2\sqrt{h_1 h_2}), \qquad h_1, h_2 \in (0, \infty).$$

(ii) *Under $\mathbb{P}_0$, the processes $(\ell_{T_b^h}(b+x))_{x \in \mathbb{N}_0}$ and $(\ell_{T_b^h}(-x))_{x \in \mathbb{N}_0}$ are time-homogeneous discrete-time Markov chains on $[0, \infty)$ with transition probabilities given by*

$$P^\star(h_1, dh_2) = e^{-h_1} \delta_0(dh_2) + e^{-h_1 - h_2} \sqrt{\frac{h_1}{h_2}} I_0'(2\sqrt{h_1 h_2})\, dh_2,$$

(4.4)
$$h_1, h_2 \in [0, \infty).$$

(iii) *The three Markov chains in* (i) *and* (ii) *are independent.*

We note that Theorem 4.1(ii) and an outline of its proof can be found in [30], (3.1-2). This proof uses an embedding of the random walk into a Brownian motion and the Brownian Ray–Knight theorem; we expect that Theorem 4.1(i) and (iii) can also be proved along these lines. In the extended version [7], using the density formula of Theorem 2.1, we provide a proof of Theorem 4.1 that is independent of the Brownian Ray–Knight theorem. This opens up the possibility of producing a new proof of this theorem, via a diffusion approximation of the Markov chains having the transition densities in (4.3) and (4.4). Furthermore, we emphasize that our proof can also be adapted to continous-time random walks on cycle-free graphs and has some potential to be extended to more general graphs. Theorem 2.1 contains far-ranging generalizations of the Ray–Knight idea, which are to be studied in future.

**Acknowledgments.** W. König would like to thank the German Science Foundation for awarding a Heisenberg grant (realized in 2003–2004). This project was initiated during an extensive visit of R. van der Hofstad to the University of British Columbia, Vancouver, Canada.

D. Brydges
Department of Mathematics
University of British Columbia
Vancouver, British Columbia
Canada V6T 1Z2
URL: http://www.math.ubc.ca/~db5d/

R. van der Hofstad
Department of Mathematics
 and Computer Science
Technical University Eindhoven
PO Box 513
5600 MB Eindhoven
The Netherlands
E-mail: rhofstad@win.tue.nl
URL: http://www.win.tue.nl/~rhofstad/

W. König
Mathematisches Institut
Universität Leipzig
Augustusplatz 10/11
D-04109 Leipzig
Germany
E-mail: koenig@math.uni-leipzig.de
URL: http://www.math.uni-leipzig.de/~koenig/